\documentclass{amsproc}
\usepackage{amssymb}
\usepackage{amscd}
\usepackage{amsmath}
\newtheorem{proposition}{Proposition}
\newtheorem{remark}{Remark}
\newtheorem{theorem}{Theorem}

\begin{document}

\title[Construction of solutions of Toda lattices by the classical moment problem.]
{Construction of solutions of Toda lattices by the classical
moment problem.}

\author{Alexander Mikhaylov } 
\address{St. Petersburg   Department   of   V.A. Steklov    Institute   of   Mathematics
of   the   Russian   Academy   of   Sciences, 7, Fontanka, 191023
St. Petersburg, Russia and Saint Petersburg State University,
St.Petersburg State University, 7/9 Universitetskaya nab., St.
Petersburg, 199034 Russia.}

\email{mikhaylov@pdmi.ras.ru}

\author{Victor Mikhaylov} 
\address{St.Petersburg   Department   of   V.A.Steklov    Institute   of   Mathematics
of   the   Russian   Academy   of   Sciences, 7, Fontanka, 191023
St. Petersburg, Russia and Saint Petersburg State University,
St.Petersburg State University, 7/9 Universitetskaya nab., St.
Petersburg, 199034 Russia.} \email{ftvsm78@gmail.com}

\begin{abstract}
Making use of formulas of J. Moser for a finite-dimensional Toda
lattices, we derive the evolution law for moments of the spectral
measure of the semi-infinite Jacobi operator associated with the
Toda lattice. This allows us to construct solutions of
semi-infinite Toda lattices for a wide class of unbounded initial
data by using well-known results from the classical moment problem
theory.

\end{abstract}


\maketitle

\section{Introduction}

The semi-infinite Toda lattice can be written \cite{T,JT} as the
following infinite-dimensional nonlinear system:
\begin{equation}
\label{Toda_eq}
\begin{cases}
\dot a_n(t)=a_n(t)\left(b_{n+1}(t)-b_n(t)\right),\\
\dot b_n(t)=2\left(a_n^2(t)-a_{n-1}^2(t)\right),\quad t\geqslant
0,\, n=1,2,\ldots,
\end{cases}
\end{equation}
for which one looks for a solution satisfying the initial
condition
\begin{equation}
\label{Toda_init} a_n(0)=a_n^0,\quad b_n(0)=b_n^0,\quad
n=1,2,\ldots.
\end{equation}
where $a_n^0,$ $b_n^0$ are real and $a_n^0>0,$ $n=1,2,\ldots$.
Methods of computing of functions $a_n(t)$, $b_n(t)$ are subjects
of many investigations, see for example \cite{T,JT,FT} and
references therein. In these papers the authors used inverse
scattering technique, which imposes  essential restrictions on
initial data, most used assumption is that $a_n^0,$ $b_n^0$
$n=1,2,\ldots$ are bounded. At the same time the question on the
possibility to construct a solution to (\ref{Toda_eq}),
(\ref{Toda_init}) for "unbounded" initial data is important
\cite{IV1,IV2}. In the present paper we make an attempt to
construct a solution to (\ref{Toda_eq}), (\ref{Toda_init}) for
quite general class of unbounded initial sequences.

We introduce two operators acting in  $l^\infty$, with domains
\begin{equation*}
D(H(t))=D(P(t))=\left\{\varkappa=(\varkappa_0,\varkappa_1,\ldots)\,|\,\varkappa_n=0,\,\text{for
$n\geqslant N_0\in \mathbb{N}$}\right\},
\end{equation*}
by the rules:
\begin{align*}
&\left(H(t)f\right)_1=
a_1(t)f_{2}+b_1(t)f_1,\\
&\left(H(t)f\right)_n=a_n(t)f_{n+1}+a_{n-1}f_{n-1}+b_n(t)f_n,\quad n=2,\ldots,\\
&\left(P(t)f\right)_1=a_1(t)f_2,\\
&\left(P(t)f\right)_n=a_n(t)f_{n+1}-a_{n-1}(t)f_{n-1},\quad
n=2,\ldots.
\end{align*}
Note that the operator $H(t)$ is given by the semi-infinite Jacobi
matrix (we keep the same notation for it):
\begin{equation}
\label{Jacobi_matr}
H(t)=\begin{pmatrix} b_1(t) & a_1(t)& 0 & \cdot & 0& 0 & \cdot\\
 a_1(t) & b_2(t)& a_2(t)& \cdot& 0& 0& \cdot& \\
\cdot & \cdot & \cdot & \cdot &  0 &  0& \cdot \\
0& 0& 0& a_{N-1}(t)& b_N(t)& a_{N}(t)& \cdot\\
0& 0& 0& 0& a_N(t)& b_{N+1}(t)& \cdot\\
\cdot& \cdot& \cdot& \cdot& \cdot& \cdot& \cdot
\end{pmatrix}.
\end{equation}
It is a well known fact \cite{JT,T} that the system
(\ref{Toda_eq}) is equivalent to the following operator equation:
\begin{equation}
\label{Toda_eqviv}
\frac{dH}{dt}=PH-HP.
\end{equation}
Let $d\rho^t(\lambda)$ denotes a spectral measure of the operator
corresponding to Jacobi matrix $H(t)$ (see \cite{Ah,S} and note
that it is defined in nonunique way if $H(t)$ is in the limit
circle case). The moments of $d\rho^t(\lambda)$ are introduced by
the rule
\begin{equation}
\label{Moment_eq} s_k(t)=\int_{-\infty}^\infty
\lambda^k\,d\rho^t(\lambda),\quad k=0,1,2,\ldots
\end{equation}

We briefly outline our motivation to treat the moments of spectral
measure as "inverse data" and study the evolution of moments with
respect to the parameter $t$. In \cite{MM6,MM7,MM8,MM9,MM10} the
authors studied forward and inverse dynamic problems for a
dynamical system with discrete time associated with finite and
infinite Jacobi matrices (we omit parameter $t$ here). For real
$b_k$, $k=1,2,\ldots$ and $a_k>0,$ $k=0,1,\ldots$ one considers
\begin{equation}
\label{Jacobi_dyn} \left\{
\begin{array}l
u_{n,\,t+1}+u_{n,\,t-1}-a_{n}u_{n+1,\,t}-a_{n-1}u_{n-1,\,t}-b_nu_{n,\,t}=0,\quad n\in \mathbb{N},\, t\in \mathbb{N}_0,\\
u_{n,\,-1}=u_{n,\,0}=0,\quad n\in \mathbb{N}, \\
u_{0,\,t}=f_t,\quad t\in \mathbb{N}_0.
\end{array}\right.
\end{equation}
This system is a discrete analog of an initial boundary value
problem for a wave equation with a potential on a half-line with
the Dirichlet control at $n=0$. The solution to (\ref{Jacobi_dyn})
is denoted by $u^f_{n,\,t}$. We fix some positive integer $T$ and
denote by $\mathcal{F}^T$ the \emph{outer space} of the system
(\ref{Jacobi_dyn}), the space of controls (inputs):
$\mathcal{F}^T:=\mathbb{R}^T$, $f\in \mathcal{F}^T$,
$f=(f_0,\ldots,f_{T-1})$. The input $\longmapsto$ output
correspondence in the system (\ref{Jacobi_dyn}) is realized by a
\emph{response operator}: $R^T:\mathcal{F}^T\mapsto \mathbb{R}^T$
defined by the rule
\begin{equation*}
\left(R^Tf\right)_t=u^f_{1,\,t}, \quad t=1,\ldots,T.
\end{equation*}
This operator (the discrete version of a dynamic
Dirichlet-to-Neumann map) plays the role of inverse data, and in
\cite{MM12,MM8} several methods of recovering the matrix $H$ from
this operator was proposed. For
$f=(f_0,f_1,\ldots),\,g=(g_0,g_1,\ldots)$ the convolution
$c=f*g=(c_0,c_1,\ldots)$ is defined by the formula
\begin{equation*}
c_t=\sum_{s=0}^{t}f_sg_{t-s},\quad t\in \mathbb{N}\cup \{0\}.
\end{equation*}
It has been shown that the response operator has the form
\begin{equation*}
\left(R^Tf\right)_t= r*f_{\cdot-1},
\end{equation*}
where the convolution kernel of $R^T$, called a \emph{response
vector}, admits the spectral representation
\begin{equation}
r_{k-1}=\int_{-\infty}^\infty
\mathcal{T}_k(\lambda)\,d\rho(\lambda),\quad k\in
\mathbb{N},\label{Resp_spectr}
\end{equation}
where $d\rho(\lambda)$ is a spectral measure of the operator
corresponding to Jacobi matrix $H$ in (\ref{Jacobi_dyn}), and
$\mathcal{T}_l(2\lambda)$, $l=1,2,\ldots$ are Chebyshev
polynomials of the second kind, i.e. they are obtained as a
solution to the following difference Cauchy problem:
\begin{equation}
\label{Chebysh} \left\{
\begin{array}l
\mathcal{T}_{t+1}+\mathcal{T}_{t-1}-\lambda \mathcal{T}_{t}=0,\\
\mathcal{T}_{0}=0,\,\, \mathcal{T}_1=1.
\end{array}
\right.
\end{equation}
Let for $K\in \mathbb{N}$ the matrix $\Lambda_K\in \mathbb{M}^{K}$
be defined by the following rule:
\begin{equation}
\label{LambdaMatr}
\Lambda_K=a_{ij}=\begin{cases} 0,\quad \text{if $i>j$},\\
0,\quad \text{if $i+j$ is odd,}\\
C_{\frac{i+j}{2}}^j(-1)^{\frac{i+j}{2}+j},
\end{cases}
\end{equation}
where $C_n^k$ are binomial coefficients. The formula
(\ref{Resp_spectr}) shows that entries of the response vector are
related to moments by the rule:
\begin{equation}
\label{Resp_and_moments}
\begin{pmatrix}
r_0\\
r_1\\
\ldots \\
r_{K-1}
\end{pmatrix}=\Lambda_K\begin{pmatrix}
s_0\\
s_1\\
\ldots \\
s_{K-1}
\end{pmatrix}.
\end{equation}

In the papers mentioned, inverse dynamic and spectral problems for
(\ref{Jacobi_dyn}) were studied, the relation of the Boundary
control method \cite{B07,B17} with the de Branges method
\cite{DBr} was established and the dynamic approach to classical
moment problems \cite{A,S} was proposed. The motivation to use the
set of moments as data, evolution of which can be used in solving
Toda system (\ref{Toda_eq}), (\ref{Toda_init}), comes from
formulas (\ref{Resp_spectr}) and (\ref{Resp_and_moments}). These
formulas say that the knowledge of response vector $r$, i.e.
dynamic inverse data for (\ref{Jacobi_dyn}) is equivalent to
knowledge of the set of moments. In other words, in view of
(\ref{Resp_spectr}), the set of moments can be treated as inverse
dynamic data for the system (\ref{Jacobi_dyn}). Main results of
the present paper concerns the evolution in time of moments
$s_k(t)$, $k=1,2,\ldots$, given by (\ref{Moment_eq}), with the
measure $d\,\rho^t(\lambda)$ being a spectral measure of $H(t)$.

In the second section we provide the necessary information on
finite-dimensional Toda lattices and adopt and rewrite some of the
results from \cite{M75} in a form, convenient for our purposes. In
the third section we remind the reader some basic facts on
classical moment problems and their relationships with Jacobi
operators and de Branges spaces. In the last section we derive
recurrent and exact formulas for the evolution  of moments
$s_k(t)$, $k=1,2,\ldots$ under the Toda flow, with the help of
which it makes it possible to construct the solution of
(\ref{Toda_eq}), (\ref{Toda_init}) for some classes of unbounded
initial data.

\section{Finite Toda lattice. Moser formula.}

We consider the initial value problem for the finite Toda lattice:
\begin{equation}
\label{Toda_eq_fin}
\begin{cases}
\dot a_{N,\,n}(t)=a_{N,\,n}(t)\left(b_{N,\,n+1}(t)-b_{N,\,n}(t)\right),\\
\dot b_{N,\,n}(t)=2\left(a_{N,\,n}^2(t)-a_{N,\,n-1}^2(t)\right),
\end{cases},\quad
t\geqslant 0,\, n=1,2,\ldots,N.
\end{equation}
where one looks for a solution satisfying the initial condition
\begin{equation}
\label{Toda_init_fin} a_{N,\,n}(0)=a_{N,\,n}^0,\quad
b_{N,\,n}(0)=b_{N,\,n}^0,\quad n=1,\ldots,N,
\end{equation}
where $a_{N,\,n}^0,$ $b_{N,\,n}^0$ are real and $a_{N,\,n}^0>0,$
$n=1,2,\ldots,N$. It is a well known fact \cite{JT,T} that the
system (\ref{Toda_eq_fin}) is equivalent to (\ref{Toda_eqviv})
with the matrix $H_N(t)$ given by $N\times N$ block of
(\ref{Jacobi_matr}) and $P_N(t)$ defined by
\begin{align*}
&\left(P_N(t)f\right)_1=a_1(t)f_2,\\
&\left(P_N(t)f\right)_n=a_n(t)f_{n+1}-a_{n-1}(t)f_{n-1},\quad n=2,\ldots,N-1, \\
&\left(P_N(t)f\right)_N=a_{N-1}(t)f_{N-1}.
\end{align*}
By $\left\{\lambda_{N,\,k}(t),\,\phi_{N,\,k}(t)\right\}_{k=1}^N$
we denote eigenvalues and eigenvectors of $H_N(t)$:
\begin{equation}
\label{H_eig}
H_N(t)\phi_{N,\,k}(t)=\lambda_{N,\,k}(t)\phi_{N,\,k}(t),\quad
\phi_{N,\,k}(t)\in \mathbb{R}^N,\quad
\left\{\phi_{N,\,k}(t)\right\}_{1}=1,
\end{equation}
and by $\rho^t_N(\lambda)$, the spectral measure of $H_N(t)$,
given by the formula
\begin{equation}
\label{Spec_mes} d\rho^t_N(\lambda)=\sum_{k=1}^N
\sigma_{N,\,k}^2(t)\delta(\lambda-\lambda_{N,\,k}(t)),
\end{equation}
where
\begin{equation*}
\sigma_{N,\,k}^2(t)=\frac{1}{\left(\phi_{N,\,k}(t),\phi_{N,\,k}(t)\right)},
\quad k=1,\ldots,N,
\end{equation*}
and $(\cdot,\cdot)$ is a scalar product in $\mathbb{R}^N$.

Below we adapt some of results from \cite{M75} (see also
\cite{MM11}) to the form convenient for us. For simplicity we
usually omit the argument $t$.
\begin{proposition}
\label{Prop_eigen} The eigenvalues of the matrix $H_N(t)$ do not
depend on $t$:
$\lambda_{N,\,j}(t)=\lambda_{N,\,j}(0)=\lambda_{N,\,j}$.
\end{proposition}
This fact follows from the representation
$\frac{dH_N}{dt}=i\left(H_NiP_N-(iP_N)H_N\right)=\{-iP_N,H_N\}$,
and thus $H_N(t)=e^{P_Nt}H_N(0)e^{-P_Nt}$.

By $\|\cdot\|$ we denote the norm in $\mathbb{R}^N$. The Weyl
function associated with $H_N$ \cite{SG,MMS} is introduced by the
rule
\begin{equation*}
m_N(\lambda):=\left(R(\lambda)e_1,e_1\right),
\end{equation*}
where
\begin{equation*}
R(\lambda)=\left(H_N(t)-\lambda I\right)^{-1},\quad
e_i=(0,\ldots,0,1,0,\ldots,0),
\end{equation*}
with $1$ being at $i-$th place.
\begin{proposition}
The following relation holds
\begin{equation}
\label{M_eqn} \frac{d}{dt}m_N(\lambda)=2a_{N,\,1}R_{21}(\lambda).
\end{equation}
\end{proposition}
\begin{proof}
We evaluate:
\begin{equation*}
\frac{dR}{dt}=-R\frac{dH}{dt}R=-RPHR+RHPR=-RP(I+\lambda
R)+(I+\lambda R)PR=PR-RP.
\end{equation*}
Using this relation we derive that
\begin{equation*}
\frac{d}{dt}m(\lambda)=\left(\left(PR-RP\right)e_1,e_1\right)=-2\left(RP
e_1,e_1\right)=2a_1R_{21}.
\end{equation*}
\end{proof}
We introduce the notation
\begin{equation*}
B_N=H_N-\lambda I=\begin{pmatrix} b_{N,\,1}-\lambda & a_{,\,1}& 0& 0& 0\\
 a_{N,\,1} & b_{N\,2}-\lambda & a_{N,\,2}& 0& 0\\
 0& a_{N,\,2} & b_{N\,3}-\lambda & a_{,\,3}& 0\\
\cdot& \cdot& \cdot& \cdot& \cdot \\
0& 0& 0& a_{N,\,N-1}& b_{N,\,N}-\lambda
\end{pmatrix}.
\end{equation*}
By $B_k$, $1\leqslant k < N$ we denote blocks of $B_N$, given by
the intersection of $k$ rows $N-k+1,\ldots,N-1,N$ and $k$ columns
$N-k+1,\ldots,N-1,N$. Introduce the notation
$\Delta_k:=\det{B_k}$, $k=1,\ldots,N$. Then using linear algebra
one can see that
\begin{align*}
m_N(\lambda)&=R_{11}=\frac{\Delta_{N-1}}{\Delta_N},\\
R_{21}(\lambda)&=R_{12}(\lambda)=\left(R(\lambda)
e_1,e_2\right)=-\frac{a_{N,\,1}\Delta_{N-2}}{\Delta_N}.
\end{align*}
These formulas allows one to rewrite (\ref{M_eqn}) in the
following form
\begin{equation}
\label{M_eqn2}
\frac{d}{dt}m_N(\lambda)=2\left(1-\left(b_{N,\,1}-\lambda\right)m_N(\lambda)\right).
\end{equation}
Representations of a Weyl function \cite{SG,MMS} and a spectral
measure (\ref{Spec_mes}) imply that
\begin{equation*}
m_N(\lambda)=\int_{R}\frac{1}{\lambda-z}\,d\rho(z)=\sum_{k=1}^N\frac{\sigma_{N,\,k}^2(t)}{\lambda-\lambda_{N,\,k}}.
\end{equation*}
Plugging the above representation into (\ref{M_eqn2}) yields the
following relation
\begin{equation*}
\sum _{k=1}^N
\frac{2\dot\sigma_{N,\,k}\sigma_{N,k}}{\lambda-\lambda_{N,\,k}}=
2\left(1-\left(b_{N,\,1}-\lambda\right)\sum_{k=1}^N\frac{\sigma_{N,\,k}^2}{\lambda-\lambda_{N,\,k}}\right),
\end{equation*}
where by dot we denote the differentiation with respect to $t$.
Multiplying the last equality by
$\left(\lambda-\lambda_{N,\,k}\right)$ and setting
$\lambda=\lambda_{N,\,k}$, we come to the following system:
\begin{equation}
\label{Sigma_syst}
\dot\sigma_{N,\,k}(t)=-(b_{N,\,1}-\lambda_{N,\,k})\sigma_{N,\,k}(t),\quad
k=1,\ldots,N.
\end{equation}
\begin{proposition}
The coefficient $b_{N,\,1}$ admits the representation
\begin{equation*}
b_{N,\,1}=\sum_{k=1}^N\lambda_{N,\,k}\sigma_{N,\,k}^2.
\end{equation*}
\end{proposition}
\begin{proof}
Denote by $C^k$, $k=1,\ldots,N$ the eigenvectors of $H_N$:
\begin{equation*}
H_NC^k=\lambda_{N,\,k}C^k,\quad C^k=\begin{pmatrix} C^k_1 \\ C^k_2
\\ \ldots \\ C^k_N\end{pmatrix},\quad k=1,\ldots,N,
\end{equation*}
such that $\|C^k\|=1$, $k=1,\ldots,N$. Then by the spectral
theorem
\begin{equation*}
C^*H_NC=\begin{pmatrix} \lambda_{N,\,1} & 0 &0 & \ldots & 0 \\
0& \lambda_{N,\,2} &0 &\ldots & 0\\
\ldots & \ldots & \ldots & \ldots & \ldots \\
 0& 0& 0& \ldots& \lambda_{N,\,N}
\end{pmatrix},\quad \text{where}\quad
C=\left(C^1|C^2|\ldots|C^N\right),
\end{equation*}
i.e. the matrix $C$ is constructed from columns $C^k$,
$k=1,\ldots,N$. Then
\begin{equation*}
H_N=C\begin{pmatrix} \lambda_{N,\,1} & 0 &0 & \ldots & 0 \\
\lambda_{N,\,2} &0 &\ldots & 0\\
\ldots & \ldots & \ldots & \ldots & \ldots \\
 0& 0& 0& \ldots& \lambda_{N,\,N}
\end{pmatrix}C^*,
\end{equation*}
from where and (\ref{Jacobi_matr}) we have that
\begin{equation*}
b_{N,\,1}=\left\{H_N\right\}_{11}=\sum_{k=1}^N
\lambda_{N,\,k}\left(C^k_1\right)^2=\sum_{k=1}^N
\lambda_{N,\,k}\left(\sigma_{N,\,k}\right)^2,
\end{equation*}
where we used the fact that $C^k=C^k_1\phi_{N,k}$ (see
(\ref{H_eig})).

\end{proof}
The above proposition allows us to rewrite the system
(\ref{Sigma_syst}) in the following form:
\begin{equation}
\label{Sigma_syst2}
\dot\sigma_{N,\,k}(t)=-\left(\sum_{j=1}^N\lambda_{N,\,j}\sigma_{N,\,j}^2(t)-\lambda_{N,\,k}\right)\sigma_{N,\,k}(t),\quad
k=1,\ldots,N.
\end{equation}
\begin{proposition}
The solution of (\ref{Sigma_syst2}) is given by the \emph{Moser
formula}:
\begin{equation}
\label{Moser_f}
\sigma_{N,\,k}^2(t)=\frac{\sigma_{N,\,k}^2(0)e^{2\lambda_{N,\,k}t}}{\sum_{j=1}^N\sigma_{N,\,j}^2(0)e^{2\lambda_{N,\,j}t}}.
\end{equation}
\end{proposition}

\section{Moments of the spectral measure of finite and semi-infinite
Jacobi operators and de Branges spaces.} \label{DB}

For a given sequence of positive numbers $\{c_1,c_2,\ldots\}$ and
real numbers $\{d_1, d_2,\ldots \}$, by $A$ we denote a
semi-infinite Jacobi matrix:
\begin{equation}
\label{Jac_matr}
A=\begin{pmatrix} d_1 & c_1 & 0 & 0 & 0 &\ldots \\
c_1 & d_2 & c_2 & 0 & 0 &\ldots \\
0 & c_2 & d_3 & c_3 & 0 & \ldots \\
\ldots &\ldots  &\ldots &\ldots & \ldots &\ldots
\end{pmatrix}.
\end{equation}
For $N\in \mathbb{N}$, by $A_N$ we denote the $N\times N$ Jacobi
matrix which is a block of (\ref{Jac_matr}) consisting of the
intersection of first $N$ columns with first $N$ rows of $A$.
Introduce the operator $A_{N}:\mathbb{R}^N\mapsto \mathbb{R}^N$ by
the rule:
\begin{equation}
\label{Jac_finite}
(A_N\psi)_n=\begin{cases}d_1\psi_1+c_1\psi_2,\quad n=1,\\
c_{n}\psi_{n+1}+c_{n-1}\psi_{n-1}+d_n\psi_n,\quad
2\leqslant n\leqslant N-1,\\
c_{N-1}\psi_{N-1}+d_N\psi_N,\quad n=N.
\end{cases}
\end{equation}
Let $d\mu_N(\lambda)$ denotes the spectral measure of $A_N$,
constructed by the formula (\ref{Spec_mes}).

With the semi-infinite matrix $A$ we associate the symmetric
operator $A$ (we keep the same notation) in the space $l_2$,
defined on finite sequences:
\begin{equation*}
D(A)=\left\{\varkappa=(\varkappa_0,\varkappa_1,\ldots)\,|\,\varkappa_n=0,\,\text{for
$n\geqslant N_0\in \mathbb{N}$}\right\},
\end{equation*}
and given by the rule
\begin{align*}
(A\theta)_1&=d_1\theta_1+c_1\theta_2,\\
(A\theta)_n&=c_{n}\theta_{n+1}+c_{n-1}\theta_{n-1}+d_n\theta_n,\quad
n\geqslant 2. 
\end{align*}
By $[\cdot,\cdot]$ we denote the scalar product in $l_2$. For a
given sequence
$\varkappa=\left(\varkappa_1,\varkappa_2,\ldots\right)$ we define
a new sequence
\begin{align*}
&\left(G\varkappa\right)_1=
d_1\varkappa_1+c_1\varkappa_2,\\
&\left(G\varkappa\right)_n=c_{n}\varkappa_{n+1}+c_{n-1}\varkappa_{n-1}+d_n\varkappa_n,\quad
n\geqslant 2.
\end{align*}
The adjoint operator $A^*\varkappa=G\varkappa$ is defined on the
domain
\begin{equation*}
D\left(A^*\right)=\left\{\varkappa=(\varkappa_0,\varkappa_1,\ldots)\in
l_2\,|\, (G\varkappa)\in l_2\right\}.
\end{equation*}
In the limit point case (i.e when $A$ has deficiency indices
$(0,0)$), $A$ is essentially self-adjoint. In the limit circle
case (i.e. when $A$ has deficiency indices $(1,1)$) we denote by
$p(\lambda)=(p_1(\lambda),p_2(\lambda),\ldots)$,
$q(\lambda)=(q_1(\lambda),q_2(\lambda),\ldots)$ two solutions of
the difference equation (we set here $c_0$=1):
\begin{equation}
\label{Phi_def}
c_n\phi_{n+1}+c_{n-1}\phi_{n-1}+d_n\phi_n=\lambda\phi_n,\quad
n\geqslant 1,
\end{equation}
satisfying Cauchy data $p_1(\lambda)=1,$
$p_2(\lambda)=\frac{\lambda-d_1}{c_1}$, $q_1(\lambda)=0,$
$q_2(\lambda)=\frac{1}{c_1}$. Then \cite[Lemma 6.22]{Schm}
\begin{equation*}
D\left(A^*\right)=D(\overline A)\dot +\mathbb{R}p(0)\dot
+\mathbb{R}q(0),
\end{equation*}
where $\dot +$ denotes the direct sum and $\overline A$ is a
closure of $A$. All self-adjoint extensions of $A$ are
parameterized by $h\in \mathbb{R}\cup\{\infty\}$, are denoted by
$A_{\infty,\,h}$ and are defined on the domain
\begin{equation*}
D(A_{\infty,\,h})=\begin{cases} D(\overline A)\dot +\mathbb{R}(q(0)+hp(0)),\quad h\in \mathbb{R}\\
D(\overline A)\dot +\mathbb{R}p(0),\quad h=\infty.
\end{cases}
\end{equation*}
All the details the reader can find in \cite{S,Schm}. We introduce
the measure
$d\mu_{\infty,\,h}(\lambda)=\left[dE^{A_{\infty,\,h}}_\lambda
e_1,e_1\right]$, where $dE^{A_{\infty,h}}_\lambda$ is the
projection-valued spectral measure of $A_{\infty,\,h}$ such that
$E^{A_{\infty,\,h}}_{\lambda-0}=E^{A_{\infty,\,h}}_{\lambda}$. The
results of \cite{A} and \cite[Section 5]{S} imply that $d\mu_N\to
d\mu_{\infty,\,\alpha}$ $*-$weakly as $N\to\infty,$ where
\begin{equation}
\label{alpha_def}
\alpha=-\lim_{n\to\infty}\frac{q_n(0)}{p_n(0)}.
\end{equation}

Assume that we are given a set of real numbers
$\left\{s_k\right\}_{k=0}^\infty$. We denote by $C_T[X]$ the set
of polynomials of the order less than $T$. Then
$\left\{s_k\right\}_{k=0}^{2T-2}$ determines on $C_T[X]$ the
bilinear form by the rule: for $F,G\in C_T[X]$,
$F(\lambda)=\sum_{n=0}^{T-1}\alpha_n \lambda^n,$
$G(\lambda)=\sum_{n=0}^{T-1}\beta_n \lambda^n,$ one defines
\begin{equation}
\label{Moment_scal} \langle
F,G\rangle=\sum_{n,m=0}^{T-1}s_{n+m}\alpha_n{\beta_m}.
\end{equation}
Thus the quadratic form (\ref{Moment_scal}) is determined by the
following Hankel matrix:
\begin{equation}
\label{S_matr}
S_T=\begin{pmatrix} s_0 & s_1 & s_2 & \ldots & s_{T-1}\\
s_1 & s_2 & \ldots& \ldots & \ldots \\
s_2 & \ldots & \ldots & \ldots & \ldots \\
\ldots & \ldots & \ldots & \ldots & s_{2T-1}\\
s_{T-1}& \ldots & \ldots & s_{2T-1} & s_{2T-2}
\end{pmatrix}
\end{equation}
In \cite{MM6,MM7,MM8,MM9,MM10} the authors studied the inverse
dynamic problems for dynamical system with discrete time
associated with finite and semi-infinite Jacobi matrices
(\ref{Jacobi_dyn}). For this system it has been proved the
statement which is equivalent to the following theorem (see also
\cite{A,S}):
\begin{theorem}
\label{Moments_char} The numbers $s_k$, $k=0,1,2,\ldots$ are
moments of some Borel measure $d\,\rho$ on $\mathbb{R}$, i.e.
\begin{equation*}
s_k=\int_{\mathbb{R}}\lambda^k\,d\rho(\lambda),\quad
k=0,1,2,\ldots
\end{equation*}
if and only if the matrices $S_N>0$ for all $N\in \mathbb{N}$.
Then the measure $d\rho$ is a spectral measure of a Jacobi
operator associated with (\ref{Jac_matr}) (not uniquely determined
when $A$ is in the limit circle case), and the block $A_N$ of $A$
(i.e. coefficients $d_1,d_2,\ldots,d_N$, $c_1,c_2,\ldots,
c_{N-1}$) can be recovered from the set of moments
$\{s_0,s_1,\ldots,s_{2N-2}\}.$

If $S_K>0$ for $K=1,\ldots, N_0$ and $\det{S_{N_0+1}}=0$, then
there exist a finite Jacobi operator $A_{N_0}$ (\ref{Jac_finite})
such that $d\,\rho(\lambda)$ is a spectral measure of this
operator with a finite support:
$\#\operatorname{supp}\left\{d\rho(\lambda)\right\}=N_0$.
\end{theorem}

In \cite{MM7,MM8,MM9,MM10} it has been shown that $C_T[X]$ is a de
Branges \cite{DBr} space associated with the system
(\ref{Jacobi_dyn}), the scalar product in which is given by
\begin{equation}
\label{Scal_DB_X} \langle
F,G\rangle=\sum_{n,m=0}^{T-1}s_{n+m}\alpha_n{\beta_m}=\int_{-\infty}^\infty
{F(\lambda)}\overline{G(\lambda)}\,d\mu(\lambda)=\int_{-\infty}^\infty
{F(\lambda)}\overline{G(\lambda)}\,d\mu_T(\lambda),
\end{equation}
where $d\mu(\lambda)$ is (any) spectral measure associated with
semi-infinite matrix $A$,  $d\mu_T(\lambda)$ is the spectral
measure associated by (\ref{Jac_finite}) with the block $A_T$; and
$F,G\in C_T[X]$, $F(\lambda)=\sum_{n=0}^{T-1}\alpha_n \lambda^n,$
$G(\lambda)=\sum_{n=0}^{T-1}\beta_n \lambda^n,$ are connected with
controls $f,g\in\mathcal{F}^T$ in (\ref{Jacobi_dyn}) by the rules
\begin{equation*}
F(\lambda)=\sum_{k=1}^T
\mathcal{T}_k(\lambda)f_{T-k},\,G(\lambda)=\sum_{k=1}^T
\mathcal{T}_k(\lambda)g_{T-k},
\end{equation*}
where $\mathcal{T}_k$ are Chebyshev polynomials of the second kind
(\ref{Chebysh}).

\section{Evolution of moments.}

\subsection{Finite-dimensional case}

Below we use the additional subscript index $N$ for moments, to
emphasize that we are in the finite dimensional case. First we
consider the evolution of moments
\begin{equation*}
s_{N,\,k}(t)=\int_{-\infty}^\infty
\lambda^k\,d\rho^t_N(\lambda),\quad k=0,1,2,\ldots
\end{equation*}
in the case of finite-dimensional system (\ref{Toda_eq_fin}),
(\ref{Toda_init_fin}), where $d\rho^t_N(\lambda)$ is a spectral
measure of $H_N(t)$ given by (\ref{Spec_mes}). We introduce the
vector function
\begin{equation}
\label{Sigma_vect}
\Theta_N(t)=\begin{pmatrix} \widetilde \sigma_{N,\,1}(t)\\
\widetilde \sigma_{N,\,2}(t)\\ \ldots \\ \widetilde
\sigma_{N,\,N}(t)
\end{pmatrix}=\begin{pmatrix} \sigma_{N,\,1}(0) e^{\lambda_{N,\,1} t}\\
\sigma_{N,\,2}(0) e^{\lambda_{N,\,2} t}\\ \ldots \\
\sigma_{N,\,N}(0) e^{\lambda_{N,\,N} t}
\end{pmatrix}.
\end{equation}
Then (\ref{Moser_f}) and (\ref{Sigma_vect}) implies that
\begin{equation*}
\sigma_{N,\,k}(t)=\frac{\widetilde
\sigma_{N,\,k}(t)}{\|\Theta_N(t)\|},
\end{equation*}
where
\begin{equation}
\label{Sigma_til} \|
\Theta_N(t)\|=\sqrt{\sum_{j=1}^N\sigma_{N,\,j}^2(0)e^{2\lambda_{N,\,j}t}}.
\end{equation}
For $k=1,2,\ldots$ we have that
\begin{equation}
\label{SK} s_{N,\,k}(t)=\int_{\mathbb{R}}
\lambda^k\,d\rho^t_N(\lambda)=\sum_{j=1}^N
\lambda_{N,\,j}^k\sigma^2_{N,\,j}(t)=\sum_{j=1}^N
\lambda_{N,\,j}^k\frac{\widetilde\sigma^2_{N,\,j}(t)}{\|
\Theta_N(t)\|^2}.
\end{equation}
Then upon introducing the notation
\begin{equation*}
\widetilde s_{N,\,k}(t)=s_{N,\,k}(t)\|\Theta_N(t)\|^2,\quad
k=1,2,\ldots
\end{equation*}
and using (\ref{SK}) we see that for $k=1,2,\ldots$ the following
relation holds
\begin{equation}
\label{S_tilde} \dot{\widetilde s}_{N,\,k}(t)=\sum_{j=1}^N
\lambda_{N,\,j}^k
2\dot{\widetilde\sigma}_{N,\,j}(t)\widetilde\sigma_{N,\,j}(t)=\sum_{j=1}^N
\lambda_{N,\,j}^{k+1} 2\widetilde\sigma_{N,\,j}^2(t)=2\widetilde
s_{N,\,k+1}(t).
\end{equation}

By $C_N[X]$ we denote the de Branges space associated with finite
Jacobi matrix $H_N(t)$ with the scalar product defined by
(\ref{Scal_DB_X}). We take $F,G\in C_N[X]$,
$F(\lambda)=\sum_{n=0}^{N-1}\alpha_n \lambda^n,$
$G(\lambda)=\sum_{n=0}^{N-1}\beta_n \lambda^n,$ the scalar product
in $C_N([X]$  has the form:
\begin{equation*}
\langle
F,G\rangle=\sum_{n,m=0}^{N-1}s_{N,\,n+m}(t)\alpha_n{\beta_m}.
\end{equation*}
We multiply both sides of the above equality by
$\|\Theta_N(t)\|^2$ and differentiate it with respect to $t$. Then
for the right hand side we get
\begin{multline}\label{Eq1}
\left(\left[F,G\right]_{B^N(t)}\|\Theta_N(t)\|^2\right)'=\sum_{n,m=0}^{N-1}\left(\|\Theta_N(t)\|^2s_{N,\,n+m}(t)\right)'\alpha_n{\beta_m}
=\sum_{n,m=0}^{N-1}\left(\widetilde
s_{N,\,n+m}(t)\right)'\alpha_n{\beta_m}\\
=\sum_{n,m=0}^{N-1}2\widetilde
s_{N,\,n+m+1}(t)\alpha_n{\beta_m}=2\|\Theta_N(t)\|^2\sum_{n,m=0}^{N-1}
s_{N,\,n+m+1}(t)\alpha_n{\beta_m}.
\end{multline}
Diffrentiating the left hand side we arrive at
\begin{equation}
\left(\left[F,G\right]_{B^N(t)}\|\Theta_N(t)\|^2\right)'
=\left(\|\Theta_N(t)\|^2\right)'\sum_{n,m=0}^{N-1}s_{N,\,n+m}(t)\alpha_n{\beta_m}+\|\Theta_N(t)\|^2\sum_{n,m=0}^{N-1}\dot
s_{N,\,n+m}(t)\alpha_n{\beta_m}.\label{Eq2}
\end{equation}
On equating (\ref{Eq2})  and the right and side of (\ref{Eq1}), we
come to the relation
\begin{eqnarray}
\frac{\left(\|\Theta_N(t)\|^2\right)'}{\|\Theta_N(t)\|^2}\sum_{n,m=0}^{N-1}s_{N,\,n+m}(t)\alpha_n{\beta_m}+\sum_{n,m=0}^{N-1}\dot
s_{N,\,n+m}(t)\alpha_n{\beta_m}
=2\sum_{n,m=0}^{N-1} s_{N,\,n+m+1}(t)\alpha_n{\beta_m}.
\end{eqnarray}
Due to the arbitrariness of $F,G$, we can formulate the following
\begin{proposition}
The moments $s_{N,\,k}(t)$ of the measure $d\rho_N^t(\lambda)$
satisfy the following recurrent relation:
\begin{equation}
\label{Main_syst} \dot
s_{N,\,k}(t)+\left(\ln{\left\{\|\Theta_N(t)\|^2\right\}}\right)'s_{N,\,k}(t)=2s_{N,\,k+1}(t),\quad
k=0,1,\ldots.
\end{equation}
\end{proposition}
Since we know that
$s_{N,\,0}(t)=1$ for all $t$, then (\ref{Main_syst}) allows us to
determine $s_{N,\,1}(t),s_{N,\,2}(t)\ldots,$ $s_{N,\,2N-2}(t)$
recursively. Then we use the fact that the set of moments
$s_{N,\,k}(t)$, $k=0,\ldots,2N-2$ determines $N\times N$ block
$H_N(t)$ of Jacobi matrix (\ref{Jacobi_matr}) and thus
coefficients $a_{N,\,k}(t)\,,b_{N,\,k}(t)\,,a_{N,\,N}(t)$,
$k=1,\ldots,N-1$. Formulas for the reconstruction of entries of
Jacobi matrix from moments are given in \cite{Ah,MM6,MM8,S}.

\subsection{More on finite-dimensional case, infinite-dimensional case.}

By $d\rho_N^0(\lambda)$ we denote a spectral measure
(\ref{Spec_mes}) of a finite Jacobi operator $H_N(0)$
corresponding to initial data in (\ref{Toda_eq_fin}),
(\ref{Toda_init_fin}):
\begin{equation}
\label{Jacobi_matr_in}
A_{N}^0=H_N(0)=\begin{pmatrix} b_{N,\,1}^0 & a_{N,\,1}^0& 0 & \cdot & 0&\\
 a_{N,\,1}^0 & b_{N,\,2}^0& a_{N,\,2}^0& \cdot& 0& \\
\cdot & \cdot & \cdot & \cdot & \cdot  \\
0& 0& 0& a_{N,\,N-1}^0& b_{N,\,N}^0
\end{pmatrix}.
\end{equation}
We observe that the
expression for the square of (\ref{Sigma_til}) admits the
following representation
\begin{equation}
\label{Omega}
\Omega_N(t):=\|\Theta_N(t)\|^2=\sum_{j=1}^N\sigma_{N,\,j}^2(0)e^{2\lambda_{N,\,j}t}=\int_{\mathbb{R}}e^{2\lambda
t}\,d\rho_N^0(\lambda),
\end{equation}
here we used Proposition \ref{Prop_eigen}. Making use (\ref{SK}),
we obtain the following
\begin{proposition}
The moments $s_{N,\,k}(t)$, $k=0,1,\ldots$ admit the
representation
\begin{equation}
\label{SK_fim}
s_{N,\,k}(t)=\frac{\int_{\mathbb{R}}\lambda^ke^{2\lambda
t}\,d\rho_N^0(\lambda)}{\int_{\mathbb{R}}e^{2\lambda
t}\,d\rho_N^0(\lambda)},\quad k-0,1,\ldots.
\end{equation}
\end{proposition}

We introduce the matrices $L_K(t)$ by the rule
\begin{multline*}
L_K(t)=\frac{1}{\Omega_N(t)}\begin{pmatrix} e^{2tA_{N,\,0}} &  e^{2tA_{N,\,0}}A_{N,\,0}   & \ldots & e^{2tA_{N,\,0}}A_{N,\,0}^{K-1}\\
e^{2tA_{N,\,0}}A_{N,\,0} & e^{2tA_{N,\,0}}A_{N,\,0}^2 & \ldots & \ldots \\
e^{2tA_{N,\,0}}A_{N,\,0}^2 & \ldots  & \ldots & \ldots \\
\ldots & \ldots &  \ldots & e^{2tA_{N,\,0}}A_{N,\,0}^{2K-1}\\
e^{2tA_{N,\,0}}A_{N,\,0}^{N-1}&  \ldots & \ldots &
e^{2tA_{N,\,0}}A_{N,\,0}^{2K-2}
\end{pmatrix}\\
=\frac{e^{2tA_{N,\,0}}}{\Omega_N(t)}\begin{pmatrix} I &  A_{N,\,0} & A_{N,\,0}^2 & \ldots & A_{N,\,0}^{K-1}\\
A_{N,\,0} & A_{N,\,0}^2 & \ldots& \ldots & \ldots \\
A_{N,\,0}^2 & \ldots & \ldots & \ldots & \ldots \\
\ldots & \ldots & \ldots & \ldots & A_{N,\,0}^{2K-1}\\
A_{N,\,0}^{K-1}& \ldots & \ldots & \ldots & A_{N,\,0}^{2K-2}
\end{pmatrix}.
\end{multline*}
Using the spectral theorem we see that for $t\geqslant 0$
\begin{equation}
\label{Matr_fin_rel} \det{L_K(t)}>0,\quad K=0,1,\ldots,N,\quad
\text{and} \quad \det{L_{N+1}(t)}=0.
\end{equation}
The relations in (\ref{Matr_fin_rel}) and Theorem
\ref{Moments_char} imply that moments
$s_{N,\,0},s_{N,\,1}(t),\ldots, s_{N,\,2N-2}(t),\ldots$ for
$t\geqslant 0$ correspond to $N\times N$ Jacobi matrix.

Now we return to the semi-infinite problem (\ref{Toda_eq}),
(\ref{Toda_init}). Introduce the matrix corresponding to initial
data (\ref{Toda_init}):
\begin{equation}
\label{Jacobi_matr_init}
A^0=\begin{pmatrix} b_1^0 & a_1^0& 0 & \cdot & 0& 0 & \cdot\\
 a_1^0 & b_2^0& a_2^0& \cdot& 0& 0& \cdot& \\
\cdot & \cdot & \cdot & \cdot &  0 &  0& \cdot \\
0& 0& 0& a_{N-1}^0& b_N^0& a_{N}^0& \cdot\\
0& 0& 0& 0& a_N^0& b_{N+1}^0& \cdot\\
\cdot& \cdot& \cdot& \cdot& \cdot& \cdot& \cdot
\end{pmatrix}.
\end{equation}
By $A^0_{N}$ we denote the $N\times N$ block of $A^0$, the
spectral measure of the operator corresponding to $A^0_{N}$ (see
(\ref{Jac_finite})) is denoted by $d\,\rho_{N}^0(\lambda)$ and is
defined by the formula (\ref{Spec_mes}). We know that
$d\rho^0_N\to d\rho^0_{\infty,\,\alpha}$ $*-$weakly as
$N\to\infty,$ where $d\rho^0_{\infty,\,\alpha}$ is a spectral
measure of $A^0_{\infty,\,\alpha}$ with $\alpha$ given by
(\ref{alpha_def}).

Solving the problem (\ref{Toda_eq_fin}) with initial conditions
(\ref{Toda_init_fin}) given by the matrix $A^0_{N}$, we get
formulas (\ref{Omega}) for norming coefficients and (\ref{SK_fim})
for moments of $d\rho^0_N(\lambda)$.

One can observe that if in (\ref{Omega}), (\ref{SK_fim}) the
support of the "limit measure"
$d\,\rho_{\infty,\,\alpha}^0(\lambda)$ is semibounded, then it is
possible to go to the limit as $N\to\infty$. Specifically, the
following proposition holds:
\begin{proposition}
If the measure $d\rho_{\infty,\,\alpha}^0(\lambda)$ is such that
\begin{equation}
\label{measure_restr}
\operatorname{supp}\left\{d\,\rho_{\infty,\,\alpha}^0(\lambda)\right\}\subset
(-\infty,M)
\end{equation}
for some $M\in \mathbb{R}$, then there exist the following limits
\begin{align}
s_k(t)&:=\lim_{N\to\infty}s_{N,\,k}(t)=\frac{\int_{\mathbb{R}}\lambda^ke^{2\lambda
t}\,d\rho_{\infty,\,\alpha}^0(\lambda)}{\int_{\mathbb{R}}e^{2\lambda
t}\,d\rho_{\infty,\,\alpha}^0(\lambda)},\quad k=0,1,\ldots,\label{SK_def}\\
\Omega_\alpha(t)&:=\lim_{N\to\infty}\Omega_N(t)=\int_{\mathbb{R}}e^{2\lambda
t}\,d\rho_{\infty,\,\alpha}^0(\lambda).\notag
\end{align}
Moreover, the functions $s_k(t)$ satisfy the recurrent relation
\begin{equation}
\label{Main_eq} \dot
s_{k}(t)+\left(\ln{\left\{\Omega_\alpha(t)\right\}}\right)'s_{k}(t)=2s_{k+1}(t),\quad
k=0,\ldots,\,\,t>0,
\end{equation}
where $s_0(t)=1$ for $t\geqslant 0.$
\end{proposition}

Now we can treat the functions constructed by (\ref{SK_def}) (or
by (\ref{Main_eq})) as moments of a spectral measure of some
Jacobi matrix, whose coefficients (depending on $t$) we call the
solution to (\ref{Toda_eq}), (\ref{Toda_init}).
\begin{remark}
The fact that numbers $s_k(t)$ for all $t>0$ defined by
(\ref{SK_def}) are indeed moments of some Borel measure on
$\mathbb{R}$ and thus yield the Jacobi matrix (which depends on
$t$), follows from the observation that for any $t>0$,
$s_k(t)=\frac{1}{\Omega_\alpha(t)}\gamma_k(t)$, where
$\gamma_k(t)$ are the moments of the measure $e^{2\lambda
t}d\rho_{\infty,\,\alpha}^0(\lambda)$:
$\gamma_k(t)=\int_{-\infty}^\infty \lambda^k e^{2\lambda
t}\,d\rho_{\infty,\,\alpha}^0(\lambda)$. And as such, they satisfy
conditions of Theorem \ref{Moments_char}. One can also use
arguments of finite-dimensional case and Theorem
\ref{Moments_char} to show that corresponding determinants are all
positive.
\end{remark}

We conclude the paper with the following
\begin{theorem}
If the spectral measure $d\,\rho_{\infty,\,\alpha}^0(\lambda)$ of
the operator $A^0$ corresponding to Cauchy data (\ref{Toda_init})
satisfy the restriction (\ref{measure_restr}), then the solution
$a_n(t),$ $b_n(t)$, $n=1,2,\ldots$ to (\ref{Toda_eq}),
(\ref{Toda_init}) is determined by the moments $s_k(t)$,
$k=0,1,\ldots$ satisfying the recurrent relation (\ref{Main_eq})
and given by formulas (\ref{SK_def}).
\end{theorem}
Formulas for the reconstruction of $a_k(t)$, $b_k(t)$ from
$s_k(t)$ are standard and are provided in \cite{Ah,MM6,MM8,S}

\noindent{\bf Acknowledgments}

A. S. Mikhaylov and V. S. Mikhaylov were partly supported by RFBR
18-01-00269, and by Volkswagen Foundation project "From Modeling
and Analysis to Approximation".

\end{document}